\documentclass{article}

\usepackage{graphics}
\usepackage{amssymb,amsmath,array,epsfig}

\newtheorem{theorem}{Theorem}[section]

\newtheorem{e-proposition}[theorem]{Proposition}

\newtheorem{e-definition}[theorem]{Definition\rm}

\newcommand\Sym{{\mathsf{S}}}
\newcommand\Stab{{\mathsf{St}}}
\newcommand\tree{{\mathcal{T}}}

\newcommand\Ht{{H^{(3)}}}
\newcommand\Hf{{H^{(4)}}}
\newcommand\Hk{{H^{(k)}}}

\title{Asymptotic aspects of Schreier graphs and Hanoi Towers groups}

\author{Rostislav Grigorchuk~\thanks{Partially supported by NSF grants DMS-0308985 and
DMS-0456185.}\hspace{3mm}and Zoran {\v S}uni{\'k} \vspace{2mm}\\
Department of Mathematics,\\ Texas A\&M University, MS-3368, \\
College Station, TX, 77843-3368, USA }

\date{}

\begin{document}
\maketitle

\begin{abstract}
We present relations between growth, growth of diameters and the
rate of vanishing of the spectral gap in Schreier graphs of
automaton groups. In particular, we introduce a series of examples,
called Hanoi Towers groups since they model the well known Hanoi
Towers Problem, that illustrate some of the possible types of
behavior.
\end{abstract}

\section{Actions on rooted trees, automaton groups, and Hanoi Towers groups}

The free monoid $X^*$ of words over the alphabet $X=\{0,\dots,k-1\}$
ordered by the prefix relation has a $k$-regular rooted tree
structure in which the empty word is the root and the words of
length $n$ constitute the level $n$ in the tree. The $k$ children of
the vertex $u$ are the vertices $ux$, for $x=0,\dots,k-1$. Denote
this $k$-regular rooted tree by $\tree$. Any automorphism $g$ of
$\tree$ can be (uniquely) decomposed as $g = \pi_g \
(g_0,g_1,\dots,g_{k-1})$, where $\pi_g \in \Sym_k$ is called the
\emph{root permutation} of $g$, and $g_x$, $x=0,\dots,k-1$, are tree
automorphisms, called the (first level) \emph{sections} of $g$. The
root permutation $\pi_g$ and the sections $g_i$ are determined
uniquely by the relation $g(xw) = \pi_g(x) g_x(w)$, for all $x \in
X$ and $w \in X^*$. The action of a tree automorphism can be
extended to an isometric action on the boundary $\partial \tree$
consisting of the infinite words over $X$. The space $\partial
\tree$ is a compact ultrametric space homeomorphic to a Cantor set.

For any permutation $\pi$ in $\Sym_k$ define a $k$-ary tree
automorphism $a=a_\pi$ by $a = \pi \ (a_0,a_1,\dots,a_{k-1})$, where
$a_i$ is the identity automorphism if $i$ is in the support of $\pi$
and $a_i=a$ otherwise. The action of the automorphism $a_{(ij)}$ on
$\tree$ is given (recursively) by
\begin{equation}\label{aij}
a_{(ij)}(iw) = jw, \qquad\qquad a_{(ij)}(jw) = iw, \qquad\qquad
    a_{(ij)}(xw) = xa_{(ij)}(w), \text{ for }x \not\in \{i,j\}.
\end{equation}
\emph{Hanoi Towers group on $k$ pegs}, $k \geq 3$, is the group $\Hk
= \langle \{ a_{(ij)} \mid 0\leq i < j \leq k-1 \} \rangle$ of
$k$-ary tree automorphisms generated by the automorphisms
$a_{(ij)}$, $0\leq i < j \leq k-1$, corresponding to the
transpositions in $\Sym_k$. The group $\Hk$ derives its name from
the fact that it models the well known Hanoi Towers Problem on $k$
pegs (see~\cite{bode-h:open,klavzhar-m-p:frame}). In this Problem,
$n$ disks of distinct sizes, enumerated $1,\dots,n$ by their size,
are placed on $k$ pegs, denoted $0,\dots,k-1$. Any configuration of
disks is allowed as long as no disk is placed on top of a smaller
disk. A legal move consists of moving the top disk from one peg to
the top of another peg (as long as the new configuration is
allowed). The $n$-disk configurations are in bijective
correspondence with the $k^n$ words of length $n$ over
$X=\{0,\dots,k-1\}$. Namely, the word $x_1 \dots x_n$ over $X$
corresponds to the unique configuration in which the disk $i$,
$i=1,\dots,n$, is placed on peg $x_i$. The action of the
automorphism $a_{(ij)}$ corresponds to a move between the pegs $i$
and $j$.

The group $\Hk$, $k \geq 3$, is an example of an automaton group. In
general, an \emph{invertible automaton} is a quadruple
$A=(S,X,\tau,\rho)$ in which $S$ is a finite set of states, $X$ a
finite alphabet, $\tau:S \times X \to S$ a \emph{transition
function} and $\pi: S \times X \to X$ an \emph{output function} such
that, for each state $s \in S$, the restriction $\pi_s=\pi(s,.):X
\to X$ is a permutation in $\Sym_X$
(see~\cite{grigorchuk-n-s:automata}). The states of $A$ define
recursively tree automorphisms by setting the permutation $\pi_s$ to
be the root permutation of $s$ and the state $\tau(s,x)$ to be the
section $s_x$ of $s$ at $x$. The group of tree automorphisms $G(A) =
\langle s \mid s \in S \rangle$ generated by the automorphisms
corresponding to the states of the invertible automaton $A$ is
called the \emph{automaton group} of $A$. Invertible automata are
often represented by diagrams such as the one on the left in
Figure~\ref{figure} corresponding to $\Hf$. Each state $s$ is
labeled by the permutation $\pi_s$ and the labeled edges describe
the transition function (if $\tau(s,x)=t$ then there is an edge
labeled $x$ connecting $s$ to $t$).

An automaton group is \emph{contracting} if the length (with respect
to the generating set $S$) of each section $g_i$, $i=0,\dots,k-1$,
of $g$ is shorter than the length of $g$, for all sufficiently long
elements $g$. A spherically transitive (transitive on all levels)
group of tree automorphisms $G$ is \emph{regularly branching} over
its subgroup $K$ if $K$ is a normal subgroup of finite index in $G$
such that $K \times \dots \times K$ ($k$ factors) is a normal
subgroup of finite index in $G$ that is geometrically contained in
$K$ (meaning that the $k$ factors act independently on the
corresponding $k$ subtrees of $\tree$ rooted at the vertices of the
first level of $\tree$; see~\cite{bartholdi-g-s:branch}).

\begin{theorem}
The Hanoi Towers group $\Ht$ is contracting and regularly branching
over its commutator.
\end{theorem}

The groups $\Hk$ act spherically transitively on $\tree$, but are
not contracting for $k \geq 4$.

\section{Schreier graphs}

Let $A$ be an automaton, $G=G(A)$ the corresponding automaton group,
$\xi$ an infinite geodesic ray starting at the root of $\tree$,
$\xi_n$ the prefix of $\xi$ of length $n$, $P_\xi$ the stabilizer
$\Stab_G(\xi)$ and $P_n$ the stabilizer $\Stab_G(\xi_n)$,
$n=0,1,\dots$~. Denote by $\Gamma_\xi$ (or just $\Gamma$) the
Schreier graph $\Gamma = \Gamma(G,P_\xi,S)$ and by $\Gamma_n$ the
Schreier graph $\Gamma_n = \Gamma_n(G,P_n,S)$, where $S$ is the
generating set defined by the states of $A=(S,X,\tau,\rho)$. We
assume that the action of $G$ on $\tree$ is spherically transitive.
Thus the graphs $\Gamma_n$ are connected, have size $k^n$,
$n=0,1,\dots$, and $G$ is infinite. The graph $\Gamma_n$ is indeed
the graph of the action of $G$ on level $n$ in the tree and $\Gamma$
is the graph of the action on the orbit of $\xi$ in $\partial
\tree$. The sequence of graphs $\{\Gamma_n\}$ converges to the
infinite graph $\Gamma$ in the space of pointed
graphs~\cite{grigorchuk-z:cortona} (based at $\xi_n$, $n=0,1,\dots$,
and $\xi$, respectively).

\begin{figure}
\begin{center}
\begin{tabular}{ccc}
\epsfig{file=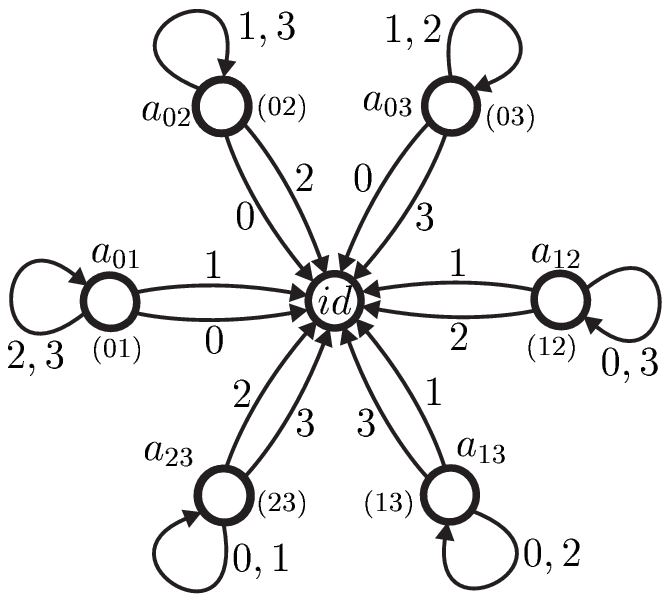,height=150pt} & \hspace{1cm}
\epsfig{file=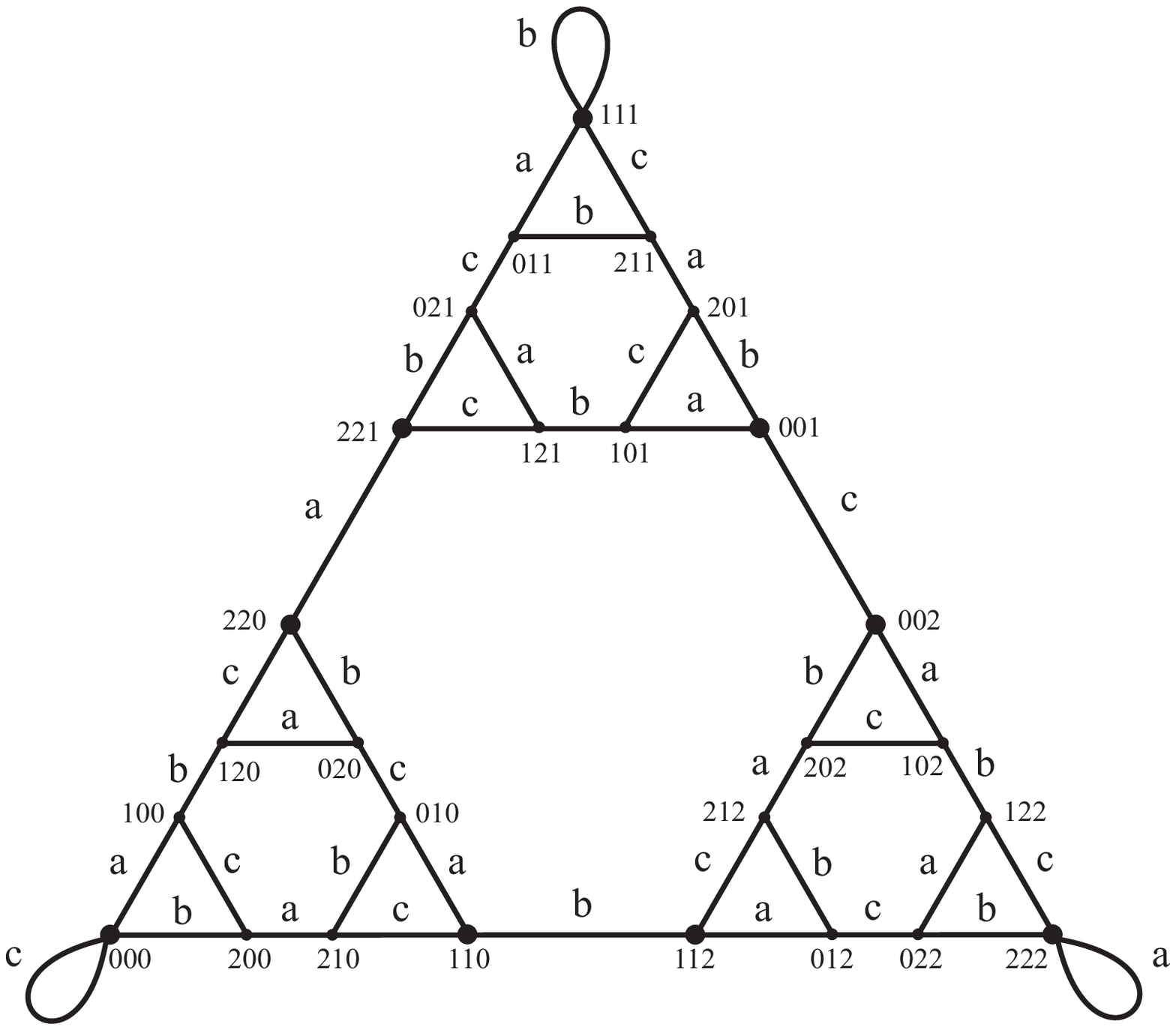,height=210pt}
\end{tabular}
\caption{The automaton generating $\Hf$ and the Schreier graph of
$\Ht$ at level 3 / L'automate engendrant $\Hf$ et le graphe de
Schreier de $\Ht$ au niveau 3}\label{figure}
\end{center}
\end{figure}

For example, in the case of the Hanoi group on 3 pegs $\Ht$ the
Schreier graph $\Gamma_3$ is given in Figure~\ref{figure} (the
automorphisms $a_{(01)}$, $a_{(02)}$ and $a_{(12)}$ are denoted by
$a$, $b$ and $c$, respectively in the figure). The group $\Ht$ has
been independently constructed in~\cite{nekrashevych:b-selfsimilar}
as a group whose limit space is the Sierpi{\'n}ski gasket.

All our asymptotic considerations are intended in the following
sense. For two functions, we write $f \preceq g$ if there exists
$C>0$ such that $f(n) \leq g(Cn)$, for all $n \geq 0$, and $f \sim
g$ if $f \preceq g$ and $g \preceq f$.

\subsection{Growth}

Denote by $\gamma(n)$ the growth function of $\Gamma$ counting the
number of vertices in $\Gamma$ that are in the ball of radius $n$
centered at $\xi$. The growth of $\Gamma$ is either exponential
($\underset{n \to \infty}{\lim} \sqrt[n]{\gamma(n)} = c>1$),
intermediate ($\underset{n \to \infty}{\lim} \sqrt[n]{\gamma(n)} =
1$, $n^m \preceq \gamma(n)$, for all $m>0$) or polynomial
($\underset{n \to \infty}{\lim} \sqrt[n]{\gamma(n)} = 1$, $\gamma(n)
\preceq n^m$, for some $m>0$). The growth is polynomial in any
Schreier graph $\Gamma$ associated to a contracting
group~\cite{bartholdi-g:spectrum}. For example, $\gamma(n) \sim
n^{\log_2(3)}$ in case of $\Ht$, and the growth is linear for the
first group from~\cite{grigorchuk:burnside} and the Erschler
group~\cite{erschler:subexponential}. The growth of $\Gamma$ is
exponential for the automata generating the lamplighter
group~\cite{grigorchuk-z:l2} and Baumslag-Solitar solvable
groups~\cite{bartholdi-s:bs}, the Bellaterra
automaton~\cite{nekrashevych:b-selfsimilar} generating the free
product of three cyclic groups of order 2, and the Aleshin
automaton~\cite{aleshin:free,vorobets:free} generating the free
group of rank 3.

\begin{theorem}
For the Hanoi Towers group $\Hk$, $k \geq 4$, the growth of
$\Gamma_{000\dots}$ is intermediate. Moreover, $a^{(\log n)^{k-2}}
\preceq \gamma(n) \preceq b^{(\log n)^{k-2}}$, for some constants
$b>a>1$.
\end{theorem}

Thus the Schreier graph $\Gamma$ is amenable, for all Hanoi Towers
groups (another example providing Schreier graphs of intermediate
growth is the group studied
in~\cite{benjamini-h:intermediate,bondarenko-c-n:intermediate}). It
is not known if, for $k \geq 4$, $\Hk$ is amenable or if it contains
free subgroups of rank 2 ($\Ht$ is amenable but not elementary
amenable).

\subsection{Diameters}

Denote by $d(n)$ the diameter of the graph $\Gamma_n$. The growth of
the diameter function $d(n)$ is either exponential, intermediate or
polynomial. Since we assumed spherical transitivity of the action,
we have that the Schreier graphs $\Gamma_n$ are connected and
$\gamma(d(n)) \geq k^n$. This implies that if either $d(n)$ or
$\gamma(n)$ grows polynomially then the other one grows
exponentially. Thus, contracting groups provide examples with
exponential diameter growth. As a concrete example, it is well known
that $d(n)=2^n-1$ for $\Ht$. It can be shown that examples of
polynomial (even linear) diameter growth are given by the the
realizations of the lamplighter group in~\cite{grigorchuk-z:l2} and
by Baumslag-Solitar solvable groups $BS(1,n)$, $n \neq \pm 1$,
in~\cite{bartholdi-s:bs}.

\begin{theorem}
For the Hanoi Towers group $\Hk$, $k \geq 3$, the diameter $d(n)$ is
asymptotically $e^{n^{\frac{1}{k-2}}}$.
\end{theorem}

Thus the Hanoi Towers groups provide examples with intermediate
diameter growth (for $k \geq 4$). The classical Hanoi Towers Problem
asks for the number of steps needed to reach the configuration $1^n$
from $0^n$. It can be shown that this distance is asymptotically
equal (it is not equal!) to the diameter $d(n)$. The Frame-Stewart
algorithm (see~\cite{klavzhar-m-p:frame}) solves the Hanoi Towers
Problem in $\sim e^{n^{\frac{1}{k-2}}}$ steps and it follows from
the work of Szegedy~\cite{szegedy:hanoi} that this is asymptotically
optimal.

\subsection{Spectra and spectral gap}

Let $T_n$ be the adjacency matrix of $\Gamma_n$ and $M_n =
\frac{1}{|S \cup S^{-1}|} T_n$ be the corresponding Markov operator
(related to the simple random walk on $\Gamma_n$). Similarly, let
$T$ be the adjacency matrix of $\Gamma$ and $M$ the corresponding
Markov operator. The \emph{spectrum} of $\Gamma_n$ (or $\Gamma$) is,
by definition, the spectrum of $M_n$ (or $M$). The regularity of
$\Gamma_n$ implies that 1 is in the spectrum of $\Gamma_n$. Let
$\delta(n)$ be the \emph{spectral gap} $1-\lambda(n)$, where
$\lambda(n)$ is the largest eigenvalue of $\Gamma_n$ different from
1. The family $\{\Gamma_n\}$ is, by definition, a family of
\emph{expanders} if there is a constant $c>0$ such that $\delta(n)
\geq c$, for all $n$. If $\{\Gamma_n\}$ do not form a family of
expanders then the spectral gap tends to 0 at a rate that is
exponential ($\underset{n \to \infty}{\lim}
\sqrt[n]{\delta(n)}=c<1$), intermediate ($\underset{n \to
\infty}{\lim} \sqrt[n]{\delta(n)}= 1$, $\delta(n) \preceq n^{-m}$,
for all $m>0$) or polynomial ($\underset{n \to \infty}{\lim}
\sqrt[n]{\delta(n)}= 1$, $n^{-m} \preceq \delta(n)$, for some
$m>0$). Chung inequality tells us that, for a graph $\Gamma=(V,E)$,
the diameter $d \leq \frac{\ln(|V|-1)}{-\ln(\lambda)}+1$, where
$\lambda$ is the largest absolute value of an eigenvalue of $\Gamma$
different from 1. In our situation $|V|=k^n$ and by using spectrum
shifting (achieved by adding loops) we obtain $d(n) \preceq
\frac{n}{\delta(n)}$. In case of an expander family $d(n) \preceq
Cn$, for some constant $C>0$, i.e.~the diameter growth is linear.
The converse is not true (the lamplighter and the Baumslag-Solitar
examples have linear diameter growth). Examples of automaton groups
producing either expanders or intermediate decay of the spectral gap
are not known yet. Examples with exponential decay are provided by
the first group in~\cite{grigorchuk:burnside} and $\Ht$. More
generally, for contracting groups the spectral gap decays
exponentially since the diameters grow exponentially. An example
with polynomial decay is provided by the lamplighter
group~\cite{grigorchuk-z:l2}.

The rate of convergence of $\delta_n$ to 0 is of special interest
and is related to the behavior of the KNS spectral measure of
$\Gamma$ and growth of F{\o}lner sets. The KNS measure $\nu$ is
limit of counting measures $\nu_n$ on
$\Gamma_n$~\cite{bartholdi-g:spectrum} ($\nu_n(B)$ is the ratio
$m(B)/k^n$, where $m(B)$ counts the eigenvalues of $\Gamma_n$ in
$B$).

We present here a full description of the spectrum of each 3-regular
Schreier graph $\Gamma_n$ modeling the Hanoi Towers Problem on 3
pegs as well as the spectrum of the 3-regular limit graph
$\Gamma_{000\dots}$ and the associated KNS spectral measure.

\begin{theorem}
Let $G$ be the Hanoi Towers group on three pegs $\Ht$. The spectrum
of $\Gamma_n$, as a set, has $3 \cdot 2^{n-1}-1$ elements and is
equal (after re-scaling by 3) to
\begin{equation}\label{specn}
 \{3\} \ \cup \ \bigcup_{i=0}^{n-1} f^{-i}(0) \ \cup \ \bigcup_{j=0}^{n-2} f^{-j}(-2),
\end{equation}
where $f$ is the polynomial $f(x) = x^2 - x -3$. The multiplicity of
the $2^i$ eigenvalues in $f^{-i}(0)$, $i=0,\dots,n-1$ is $a_{n-i}$
and the multiplicity of the $2^j$ eigenvalues in $f^{-j}(-2)$,
$j=0,\dots,n-2$ is $b_{n-j}$, where $a_i = \frac{3^{i-1}+3}{2}$ and
$b_j = \frac{3^{j-1}-1}{2}$, for $i,j \geq 1$.

The spectrum of $\Gamma_{000\dots}$, as a set, is equal (after
re-scaling by 3) to
\begin{equation}\label{spec}
 \overline{ \{3\} \ \cup \ \bigcup_{i=0}^\infty f^{-i}\{0,-2\} } =
 \overline{ \bigcup_{i=0}^\infty f^{-i}\{0\} }.
\end{equation}
It consists of a set of isolated points $I= \bigcup_{i=0}^\infty
f^{-i}\{0\}$ and its set of accumulation points $J$, which is the
Julia set of the polynomial $f$ and is a Cantor set. The KNS
spectral measure is discrete, concentrated on $\bigcup_{i=0}^\infty
f^{-i}\{0,-2\}$ and the measure of the eigenvalues in
$f^{-i}\{0,-2\}$ is $\frac{1}{6 \cdot 3^i}$, $i=0,1,\dots$~.
\end{theorem}

The method used to calculate the above spectra is a continuation of
the direction set in~\cite{bartholdi-g:spectrum}. At the core of
this method lies a renormalization principle. The action of $\Ht$ on
level $n$ of $\tree$ induces $3^n$-dimensional permutation matrix
representation $\rho_n$. The spectra of $\Gamma_n$, i.e.~of the
adjacency matrices $T_n = \rho_n(a)+ \rho_n(b) + \rho_n(c)$, are
then calculated by considering two-dimensional pencils
$\Delta_n(x,y) = T_n - xI_n + (y-1)K_n$ of matrices (where $I_n$ is
the identity matrix and $K_n$ is additional matrix aiding the
computation). The zeros of $\Delta_n(x,y)=0$ are related for various
$n$ by using iterations of a two-dimensional rational map $F:
{\mathbb R}^2 \to {\mathbb R}^2$. The spectrum of
$\Gamma_{000\dots}$ is the intersection of an invariant set of $F$
(consisting of a family of hyperbolae) with the line $y=1$. The
rational map $F$ is semi-conjugate to the one-dimensional polynomial
map $f(x)=x^2-x-3$ above, which makes possible the complete
calculation of the spectra.

A similar result is obtained (by using a different method that also
has the renormalization spirit in its root)
in~\cite{teplyaev:gasket} for the spectrum of the infinite
$4$-regular (with single exception) graph approximating the
Sierpi{\'n}ski gasket (the spectrum of this graph, as a set, was
known since the 1980's from the work of B\'ellissard). The spectrum
is obtained by iterations of a quadratic polynomial that is not
conjugate to $f$. On the other hand, iterations of a polynomial
conjugate to $f$ appear in~\cite{grabner-w:sierpinski} in the
functional equation for the Green's function of the graph obtained
by gluing two copies of the graph from~\cite{teplyaev:gasket} along
a vertex.


\begin{thebibliography}{10}

\bibitem{aleshin:free}
S.~V. Aleshin.
\newblock A free group of finite automata.
\newblock {\em Vestnik Moskov. Univ. Ser. I Mat. Mekh.}, (4):12--14, 1983.

\bibitem{bartholdi-g:spectrum}
L.~Bartholdi and R.~I. Grigorchuk.
\newblock On the spectrum of {H}ecke type operators related to some fractal
  groups.
\newblock {\em Tr. Mat. Inst. Steklova}, 231(Din. Sist., Avtom. i Beskon.
  Gruppy):5--45, 2000.

\bibitem{bartholdi-g-s:branch}
Laurent Bartholdi, Rostislav~I. Grigorchuk, and Zoran
{\v{S}}uni{\'k}.
\newblock Branch groups.
\newblock In {\em Handbook of algebra, Vol. 3}, pages 989--1112. North-Holland,
  Amsterdam, 2003.

\bibitem{bartholdi-s:bs}
Laurent Bartholdi and Zoran {\v{S}}uni{\'k}.
\newblock Some solvable automaton groups.
\newblock In {\em Topological and Asymptotic Aspects of Group Theory}, volume
  394 of {\em Contemp. Math.}, pages 11--30. Amer. Math. Soc., Providence, RI,
  2006.

\bibitem{benjamini-h:intermediate}
I.~Benjamini and C.~Hoffman.
\newblock {$\omega$}-periodic graphs.
\newblock {\em Electron. J. Combin.}, 12(1): R46, 12 pp.
  (electronic), 2005.

\bibitem{bode-h:open}
J.-P.~Bode and A.~M.~Hinz.
\newblock Results and open problems on the {T}ower of {H}anoi.
\newblock {\em Congr. Numer.}, 139:113--122, 1999.

\bibitem{bondarenko-c-n:intermediate}
I.~Bondarenko, T.~Ceccherini-Silberstein, and V.~Nekrashevych.
\newblock Amenable graphs with dense holonomy and no cocompact isometry groups.
\newblock (preprint), 2006.

\bibitem{erschler:subexponential}
Anna Erschler.
\newblock Boundary behavior for groups of subexponential growth.
\newblock {\em Ann. of Math. (2)}, 160(3):1183--1210, 2004.

\bibitem{grabner-w:sierpinski}
Peter~J. Grabner and Wolfgang Woess.
\newblock Functional iterations and periodic oscillations for simple random
  walk on the {S}ierpi\'nski graph.
\newblock {\em Stochastic Process. Appl.}, 69(1):127--138, 1997.

\bibitem{grigorchuk:burnside}
R.~I. Grigorchuk.
\newblock On {B}urnside's problem on periodic groups.
\newblock {\em Funktsional. Anal. i Prilozhen.}, 14(1):53--54, 1980.

\bibitem{grigorchuk-n-s:automata}
R.~I. Grigorchuk, V.~V. Nekrashevich, and V.~I. Sushchanski{\u\i}.
\newblock Automata, dynamical systems, and groups.
\newblock {\em Tr. Mat. Inst. Steklova}, 231(Din. Sist., Avtom. i Beskon.
  Gruppy):134--214, 2000.

\bibitem{grigorchuk-z:cortona}
R.~Grigorchuk and A.~{\.Z}uk.
\newblock On the asymptotic spectrum of random walks on infinite families of
  graphs.
\newblock In {\em Random walks and discrete potential theory},
  Sympos. Math., XXXIX, pages 188--204. Cambridge Univ. Press, Cambridge, 1999.

\bibitem{grigorchuk-z:l2}
Rostislav~I. Grigorchuk and Andrzej {\.Z}uk.
\newblock The lamplighter group as a group generated by a 2-state automaton,
  and its spectrum.
\newblock {\em Geom. Dedicata}, 87(1-3):209--244, 2001.

\bibitem{klavzhar-m-p:frame}
Sandi Klav{\v{z}}ar, Uro{\v{s}} Milutinovi{\'c}, and Ciril Petr.
\newblock On the {F}rame-{S}tewart algorithm for the multi-peg {T}ower of
  {H}anoi problem.
\newblock {\em Discrete Appl. Math.}, 120(1-3):141--157, 2002.

\bibitem{nekrashevych:b-selfsimilar}
Volodymyr Nekrashevych.
\newblock {\em Self-similar groups}, volume 117 of {\em Mathematical Surveys
  and Monographs}.
\newblock American Mathematical Society, Providence, RI, 2005.

\bibitem{szegedy:hanoi}
Mario Szegedy.
\newblock In how many steps the {$k$} peg version of the {T}owers of {H}anoi
  game can be solved?
\newblock In {\em STACS 99 (Trier)}, volume 1563 of {\em Lecture Notes in
  Comput. Sci.}, pages 356--361. Springer, Berlin, 1999.

\bibitem{teplyaev:gasket}
Alexander Teplyaev.
\newblock Spectral analysis on infinite {S}ierpi\'nski gaskets.
\newblock {\em J. Funct. Anal.}, 159(2):537--567, 1998.

\bibitem{vorobets:free}
M.~Vorobets and Y.~Vorobets.
\newblock On a free group of transformations defined by an automaton.
\newblock math.GR/601231, 2006.

\end{thebibliography}

\def\cprime{$'$}

\end{document}